%Paper: q-alg/9507004
%From: GIACHETTI@fi.infn.it
%Date: Fri, 7 Jul 1995 12:02:44 +0200 (WET-DST)
%Date (revised): Sat, 18 Nov 1995 13:35:05 +0100 (MET)

%%%%%%%%%%%%%%%%%%%%%%%%%%%%%%%%%%%%%%%%%%%%%%%%%%%%%%%%%%%%%%%%%%%%
%                                                                  %
%                      PLAIN TEX FILE                              %
%                                                                  %
%%%%%%%%%%%%%%%%%%%%%%%%%%%%%%%%%%%%%%%%%%%%%%%%%%%%%%%%%%%%%%%%%%%%

\magnification 1200
\font\abs=cmr8
\font\refit=cmti8
\font\refbf=cmbx8

\font\ccc=cmcsc10

\def\uno{{\bf 1}}
\def\tens{\otimes}
\def\id{{\rm id}}
\def\gin{{{}_{inv}\Gamma}}

\def\ii#1{\item{\phantom{1}#1 .\phantom{x}}}
\def\jj#1{\item{#1 .\phantom{x}}}
\def\incl{\hookrightarrow}
\def\fic{\widehat{\varphi}}
\def\pro #1 {{\ccc (#1) Proposition.}\phantom{X}}
\def\dfn #1 {{\ccc (#1) Definition.}\phantom{X}}
\def\cor #1 {{\ccc (#1) Corollary.}\phantom{X}}
\def\lem #1 {{\ccc (#1) Lemma.}\phantom{X}}
\def\rmks #1 {{\ccc (#1) Remarks.}\phantom{X}}
\def\rmk #1 {{\ccc (#1) Remark.}\phantom{X}}
\def\thm #1 {{\ccc (#1) Theorem.}\phantom{X}}
\def\dim {{\sl Proof.}\phantom{X}}
\def\fidi{\hskip5pt \vrule height4pt width4pt depth0pt}

\def\fq{{\cal F}}
\def\dop{{{\cal D}}}
\def\ckd{{C^k({\cal D},\gin)}}
\def\ck1{{C^{k+1}({\cal D},\gin)}}

\def\uq{{{\cal U}}}
\def\gd{{}_\Gamma\hskip-1pt \delta}
\def\uqo{{\cal U}^{\rm o}}
\def\ad{{\rm ad}}
\def\Ad{{\rm Ad}}
\def\dd{\widetilde{\Delta}}
\def\sd{\widetilde{S}}
\def\rf{\rho_{\hskip-1pt\scriptscriptstyle {\cal F}}}
\def\ru{\rho_{\hskip-.0pt\scriptscriptstyle {\cal U}}}
\def\rd{\rho_{\hskip-1pt\scriptscriptstyle {\cal D}}}
\def\rdn #1{\rho_{\hskip-1pt\scriptscriptstyle {\cal D}}^{\hskip-1pt
                   \scriptscriptstyle (\hskip-0.5pt#1\hskip-0.5pt)}}

\hsize= 15 truecm
\vsize= 22 truecm
\hoffset= 0.3 truecm
\voffset= 1.3 truecm

\baselineskip= 14 pt
\footline={\hss\tenrm\folio\hss} \pageno=1

\vglue 3truecm
\centerline{\bf QUANTUM DOUBLE AND DIFFERENTIAL CALCULI.}
\bigskip
\bigskip
\centerline{{\it
F.Bonechi, R.Giachetti, R.Maciocco, E.Sorace and M.Tarlini.}}
\bigskip
\baselineskip= 10 pt
\centerline{{\abs Dipartimento di Fisica,
Universit\`a di Firenze}}
\centerline{{\abs and}}
\centerline{{\abs Sezione INFN di Firenze}}
\medskip
\centerline{{\refit e-mail: {\abs $<$}name{\abs $>$} @ fi.infn.it}}
\bigskip
\bigskip
\baselineskip= 10 pt
{\refbf Abstract.} {\abs
We show that bicovariant bimodules as defined in [1] are in one to one
correspondence with the Drinfeld quantum double representations.
We then prove that a differential calculus associated to a
bicovariant bimodule of dimension {\refit n} is connected to the
existence of a particular ({\refit n}+1)--dimensional representation
of the double. An example of
bicovariant differential calculus on the non quasitriangular
quantum group {\refit E}${}_{\abs q}$(2) is developed.
The construction is studied in terms of Hochschild
cohomology and a correspondence between differential calculi and
1-cocycles is proved. Some differences of calculi on quantum and finite
groups with respect to Lie groups are stressed.
\medskip
\noindent Preprint q-alg/9507004.
}

\bigskip
\noindent {\abs {\refit Mathematics Subject Classification}: 16W30, 17B37}

\baselineskip= 14 pt
\bigskip

\bigskip
{\bf 1.}
An approach to the differential calculus on quantum groups
was proposed in [1] some years ago. In addition to the
obvious notion of differential $d$, the fundamental algebraic
structure on which the theory was founded is that of bicovariant
bimodule: by means of such an object the properties of differential forms
are encoded and extended to the noncommutative situation.
Much work has been done in this direction ever since: however
a general treatment and a classification
of differential calculi has been constructed
only for quantum groups obtained as deformation of semisimple Lie
groups [2,3].

In this letter we prove some results connecting
the differential calculus on Hopf algebras to the Drinfeld double [4]:
in the first place we show that bicovariant bimodules are in one
to one correspondence with the Drinfeld double representations. It is
then proved that a differential calculus associated to a bicovariant
bimodule of dimension $n$ is connected to the existence of a particular
$(n+1)$ - dimensional representation of the double.
This extension leads in a standard way to the definition of a Hochschild
cohomology of the double with values in the $n$-dimensional
representation space: we prove that each differential calculus
is associated to a 1-cocycle satisfying an additional condition
with respect to the enveloping algebra component of the double. We
finally give an equivalent characterization of differential calculi
in terms of the cohomology of the algebra of functions. In this case
the additional condition is proved to become
an invariance condition with respect to a natural action and we are
able to establish a one to one correspondence of differential
calculi with invariant 1-cocycles.
The general classification of differential calculi is therefore
reduced to the study of the representations of the double and to
a cohomological problem,
which can be performed with the more usual and efficient tools.
Moreover a supply of differential calculi is obtained by
observing that the coboundary
operator maps invariant 0-cochains into invariant 1-coboundaries.

Two final remarks are in order. In the first place, the construction
is completely independent of the quasi-triangular property
of Hopf algebras and can obviously be applied
to classical groups, both Lie and discrete or finite. Secondly,
we believe that some further investigations are deserved to the
peculiar fact that all the known differential calculi on quantum
and finite groups correspond to coboundaries, at difference with the
usual Lie group case, in which no invariant coboundary exists and the
classical differential calculus is determined by a nontrivial 1-cocycle.

The plan of this letter is as follows. In the next section we give
an essential {\it resum\'ee} of the bicovariant differential calculus
suited to our purposes. In the third one the Drinfeld double is sketched.
The fourth and the fifth sections are devoted to state and prove the
results concerning representations. In the sixth section
we briefly present a new differential calculus on
one of the deformations of the 2-dimensional Euclidean group [5], which is
not quasi-triangular and for which no result was known up to present.
The final section is devoted to the cohomological analysis of the
construction.

\bigskip
{\bf 2.}
Let $(\fq ,m,\Delta,S,\epsilon)$ be the Hopf algebra of the
representative functions on a Lie or on a quantum group.
A {\it bimodule} $\Gamma$
over $\fq$ is said to be {\it left-covariant} if there is defined a
coaction $\delta_\Gamma\,
:\,\Gamma \rightarrow \fq\tens\Gamma$ with the properties
$$\delta_\Gamma(a\,\gamma)=
\Delta(a)\,\delta_\Gamma(\gamma)\,, ~~~~~~~\delta_\Gamma(\gamma\, a)=
\delta_\Gamma(\gamma)\,\Delta(a)\,,~~~~~~~
(\epsilon\tens\id)\;\delta_\Gamma(\gamma)=\gamma\,,$$
for any $a\in \fq$ and $\gamma \in \Gamma$ . Analogously we speak of {\it
right-covariant bimodule} when there is a right coaction
$\gd\,:\,\fq\rightarrow \Gamma\tens\fq$
which satisfies the same relations, with the only replacement
of $(\id\tens\epsilon)$ in the last one.

It is proved in [1] that a left-covariant bimodule $\Gamma$ is completely
characterized by a set of elements $f_{ij}$ of the dual Hopf algebra
${\cal F}^*$ which will be identified to the quantized enveloping
algebra $\uq$. These elements are required to satisfy the following
properties:
 $$\Delta (f_{ij}) = f_{ik}\tens f_{kj}\,,~~~~~\quad
\epsilon\,(f_{ij})=\delta_{ij}\,,\eqno(1)$$
where the indices $(i,j)$ take their values in an
appropriate set $I$ -- which is assumed to be finite --
and where we adopt the convention of summing over
repeated indices. Moreover we have used the same symbols $\Delta$ and
$\epsilon$ for comultiplication and counit of the Hopf algebra
$\uq$ since any possible ambiguity is removed by looking
at the elements to which they are applied.

It is then shown that $\Gamma$ is a free left module
over $\fq$ generated by $\gin=\langle\omega_i\rangle$ with
right multiplication and left coaction respectively given by
$$   \omega_i\,b = (f_{ij}\star b)\,\omega_j\,,~~~~~~~~~~
\delta_\Gamma(a\, \omega_i) = \Delta(a)\,(\uno\tens \omega_i)\ ,
\eqno(2)$$
where $f\star a=(\id\tens f)\,\Delta(a)
=\sum_{(a)}\,a_{(1)}\langle f\,,\,a_{(2)}\rangle$. Here
$\langle\,,\,\rangle$ denotes the natural duality coupling between
$\fq$ and $\uq$. Due to the second of relations
$(2)$ the elements $\omega_i$ are said to be {\it left-invariant}.
A right-covariant bimodule has the same structure of a free left
module generated by right invariant elements $\eta_i$ with the
right multiplication induced by $\eta_i\,b = (b\star f_{ij})\,\eta_j\ ,$
where $\,b\star f = (f\tens \id)\,\Delta(b)$ .

A left- and right-covariant bimodule is said to be {\it bicovariant} if
$$(\id\tens \gd)\,\delta_\Gamma = (\delta_\Gamma\tens
\id)\,\gd\ .$$
A bicovariant bimodule is characterized by $R_{ij}\in \fq$ with
$(i,j)\in I$ such that
$$\Delta (R_{ij}) = R_{ik}\tens R_{kj}\,,~~~~~~~~~~~~
\epsilon\, (R_{ij})=\delta_{ij}\eqno(3)$$
and
$$R_{ij}\,(a\star f_{ik}) = (f_{ji}\star a)\,R_{ki}\ , \eqno(4)$$
for any $a\in\fq$. For left-invariant forms the right coaction is
defined as
$$\gd(\omega_i)=\omega_j\tens R_{ji}\,.\eqno(5)$$

The bicovariance of $\Gamma$ implies that
$\Lambda^{ij}_{k\ell}=\langle f_{j\ell}\,,\,R_{ki}\rangle$
verifies the quantum Yang-Baxter equation
$\Lambda_{12}\Lambda_{13}\Lambda_{23}=\Lambda_{23}\Lambda_{13}\Lambda_{12}$
with the assignment
$(A\tens B)_{j\ell}^{ik}=A_{ij} B_{k\ell}$.
We want to stress that this
circumstance is completely independent of the quasi-triangularity
property of the Hopf algebra  $\uq$: indeed we shall show
that the appearance of the $R$-matrix is
due to the quasi-triangular property of the double.

Suppose now that for a bicovariant bimodule $\Gamma$ there exists a set
of elements $\chi_i\in \uq$ with the following two properties:
$(i)$ they generate a vector
space ${\bf g}$ (the ``{\it quantum Lie algebra}'')
closed under the adjoint action $\ad:\uq\tens
\uq\rightarrow \uq$ defined by
$\ad_X(Y)=\sum_{(X)}S(X_{(1)})\, Y\,X_{(2)}$; $(ii)$ they satisfy
$$
\Delta(\chi_i)=\chi_j\tens f_{ji}+\uno\tens \chi_i\ ,\quad\quad
\epsilon\,(\chi_i)=0\ .$$
We then define  the {\it differential} as the linear
mapping $d:\fq\rightarrow \Gamma$ given by
$$da=(\chi_i\star a)\,\omega_i\ ,\quad\quad a\in\fq\ ,$$
and we say that the couple $(\Gamma, d)$ is a {\it bicovariant first order
differential calculus} on $\fq$. If we define the vector $\chi(a)\in
\gin$ with components $[\chi(a)]_i=\langle\chi_i\,,\,a\rangle$ the
differential can be written $da=\sum\limits_{(a)} a_{(1)} \chi(a_{(2)})$.
We remark that this construction is perfectly meaningful in the classical
case: the answer is here given by the Friedrichs theorem, which selects
uniquely the Lie algebra as the vector space ${\bf g}$.

It is straightforward to prove
that the differential satisfies the Leibniz rule $d(ab)=adb+(da)b$,
and that the right ideal
${\cal J}=\{a\in {\rm ker}\,\epsilon\, |\, \langle
\chi_i\,,\,a\rangle=0\,,~\forall\chi_i\in
{\bf g}\}$
is $\ad^*$-invariant, {\it i.e.} $\ad^*\,{\cal J}
\subseteq {\cal J}\tens\fq$,
where $\ad^*(a)=\sum_{(a)}a_{(2)}\tens
S(a_{(1)}) a_{(3)}$. Moreover, according to [1], ${\cal J}$ determines
completely the bicovariant differential calculus.

\bigskip
{\bf 3.}
Consider a linear basis $\{e_A\}$ in $\fq$, its dual
basis $\{e^A\}$ in $\uq$ ( {\it i.e.} $\langle e^A\,,\,e_B\rangle =
\delta^A_B$ ) and the canonical element  $T=e_A\tens e^A\in \fq\tens\uq$
discussed in [6].
According to Drinfeld [4] we define the {\it quantum double}
$\dop$ of the Hopf algebra $\fq$ as  the unique algebra
with the following properties:
$(i)$ it is equal to $\fq\tens\uq$ as a linear
space; $(ii)$ it contains $\fq$ and $\uqo$ as Hopf
subalgebras, where $\uqo$ is $\uq$ with opposite
comultiplication; $(iii)$ it
is quasi-triangular, its ${\cal R}$ matrix being the image of the
canonical element under the natural embedding
$\fq\tens\uq\incl \dop\tens\dop$.
The quantum double can be regarded as a quantum universal enveloping
algebra in the sense and with the algebraic procedure explained
in [4].

In the following the canonical element will be called {\it universal
${\cal R}$ matrix} when considered in $\dop\tens\dop$
and {\it universal T matrix} when considered in
$\fq\tens\uq$. Denoting by $\rf$ and
$\ru$ two finite dimensional representations of $\fq$ and $\uq$,
the matrix elements
$$T_{ij}=e_A\,[\ru(e^A)]_{ij}\in\fq\ ,\quad\quad
t_{ij}=[\rf(e_A)]_{ij}\,e^A\in \uq\ ,
$$
satisfy (1) and (3) respectively:
$$\eqalign{
\Delta T_({ij})=&\,T_{ik}\tens T_{kj}\,, \hskip3.28truecm
\epsilon\,(T_{ij})=\delta_{ij}\,,\cr
\Delta(t_{ij})=&\,t_{ik}\tens t_{kj}\,, \hskip3.0truecm
\epsilon\,(t_{ij})=\delta_{ij}\,.\cr}
$$
Conversely, it is obvious that a solution of (1) and (3)
gives a representation of
$\fq$ and $\uq$ in such a way that $f_{ij}$ and $R_{ij}$
are matrix elements of $T$.

\bigskip
{\bf 4.} We now prove our first result that relates the representations of
the Drinfeld double to bicovariant bimodules.
\smallskip
\thm 6 {\it The representations $\rf$ and $\ru$ that define a
bicovariant bimodule over $\fq$ are in one to one correspondence with the
representations $\rd$ of the Drinfeld double $\dop$ by means of the
relations
$$\rf = \rd\vert_\fq \quad\quad\quad \ru = (\rd \circ
\sd^{-1})^t\vert_{\uq}\ ,
$$
where $\sd$ is the antipode of $\dop$ and $(\,)^t$ denotes
transposition.
}
\smallskip
\dim
Let us introduce the structure constants of the double
in terms of those of $\fq$:
$$\eqalign{
e_A\, e_B =&\, m_{AB}^C\, e_C \hskip3.11truecm e^A\, e^B =
\Delta^{AB}_C\, e^C\cr
\dd(e_A)=&\, \Delta_A^{BC}\,e_B\tens e_C \hskip2truecm \dd(e^A) =
m_{CB}^A\,e^B\tens e^C\cr
\sd(e_A) =&\, S_A^B\,e_B \hskip3.4truecm \sd(e^A) = (S^{-1})_B^A\,
e^B\ .\cr} \eqno(7)
$$
The relations between the elements of $\fq$ and $\uq$ are
induced by the quasi-triangu\-larity condition $\dd':=\sigma\circ\dd=
{\cal R}\,\dd {\cal R}^{-1}$, where $\sigma$ is the usual permutation
of the tensor spaces.  Explicitly:
$$\Delta_C^{AB}\,m_{BD}^E\; e_A\, e^D =
\Delta_C^{BA}\,m_{DB}^E\; e^D\, e_A\ \eqno(8)$$

Let us suppose that $\rf$ and $\ru$ define a bicovariant
bimodule, {\it i.e.} $f_{ij}=[\rf(e_A)]_{ij}\,e^A$ and $R_{ij}=
e_A\,[\ru (e^A)]_{ij}\ $. From the relation (4) it is easy
to obtain $$\Delta_B^{AC}\, m_{EC}^D\;
[\ru (e^E)]_{ij}\; [\rf (e_A)]_{ik}=
\Delta_B^{CA}\, m_{CE}^D\; [\rf (e_A)]_{ji}\;
[\ru (e^E)]_{ki}\ .\eqno(9)$$
We write now the quasi-triangularity condition on $\dop$ in a
form that will be useful for the proof:
$${\cal R}^{-1}\,\dd'(\sd^{-1}(e^D)) = \dd(\sd^{-1}(e^D))\,{\cal
R}^{-1}\ ,\eqno(10)$$
where ${\cal R}^{-1}=\sd(e_A)\tens e^A=e_A\tens \sd^{-1}(e^A)$.

Using the Hopf algebra properties and the relations of the double
as in (8), we rewrite (10) in the following form
$$\Delta_B^{AC}\, m_{EC}^D\ \sd^{-1}(e^E)\, e_A=
\Delta_B^{CA}\, m_{CE}^D\ e_A\,\sd^{-1}(e^E)\ .\eqno(11)$$
Comparing (9) and (11) it is clear that
$$\rd\vert_\fq = \rf \quad {\rm and} \quad \rd\vert_{\uq} =
(\ru\circ S^{-1})^t$$
is a representation of $\dop$.

Conversely, starting from a representation of $\dop$,
equation (9) is satisfied and $f_{ij}=[\rd(e_A)]_{ij}\, e^A$,
$R_{ij}=e_A\,[\rd (\sd^{-1}(e^A))]_{ji}$ define a
bicovariant bimodule. \fidi
\smallskip
\rmk 12 It is finally evident that the numerical
$R$-matrix $\Lambda_{k\ell}^{ij} = \langle f_{j\ell}\,,\,R_{ki}\rangle$
that appears in the theory of bicovariant bimodules
comes from these representations of the double that is
quasi-triangular by construction. Indeed using the results of Theorem
(6) we have that $\Lambda_{k\ell}^{ij}=[\rf(e_A)]_{j\ell}\,[\ru(e^A)]_{ki}=
[\rd(e_A)]_{j\ell}\,[\rd(\sd^{-1}(e^A)]_{ik}=[\sigma\circ{\cal
R}^{-1}]_{k\ell}^{ij}$. \fidi

\bigskip
{\bf 5.}
Given a bicovariant bimodule it is not always possible to construct a
differential calculus. The second result of this letter states the
connection between the existence of a bicovariant differential
calculus of dimension $n$ and the existence of a particular
representation of dimension $n+1$ of the double.
\smallskip
\thm 13 {\it Let $\Gamma$ be an $n$-dimensional bicovariant bimodule
determined by a representation $\rdn{n}$ of  $\dop$.
Suppose there exist $n$ linearly independent elements
 $\chi_i\in \uq\ ,\ i=1,\cdots,n\ ,$ such that
$\rdn{n+1}$ defined by
$$
\rdn{n+1}(e_A)=\pmatrix{
\epsilon\,(e_A)&\langle\chi_i\,,\,e_A\rangle  \cr
                             &                              \cr
      0                      & \rdn{n}(e_A)                 \cr
                             &                              \cr}\,,
{}~~~~~
\rdn{n+1}(e^A)=\pmatrix{
\epsilon\,(e^A)&     0         \cr
                             &               \cr
      0                      & \rdn{n}(e^A)  \cr
                             &               \cr}\,,\eqno(14)
$$
\smallskip\noindent
is a $(n+1)$-dimensional representation of the double. Then $(\Gamma,d)$,
where $da=(\chi_i\star a)\,\omega_i\ ,\ a\in\fq$,  defines a
bicovariant first order differential calculus.

Conversely given a
bicovariant differential calculus on $\Gamma$ the matrices {\rm (14)}
define a representation of the double.
}
\smallskip
\dim
Observe that $\rdn{n+1}|_{\uq}$ is a representation of $\uq$,
while the $\rdn{n+1}|_{\fq}$ is a representation of $\fq$ if and
only if
$$
\Delta(\chi_i)=\chi_j\tens f_{ji}+\uno\tens \chi_i\ ,\quad\quad
\epsilon\,(\chi_i)=0\ ,
\eqno(15)$$
where the $f_{ij}=[\rdn{n}(e^A)]_{ij}\,e_A$. By making explicit
the quasi-triangularity conditions of the double
on the representation $\rdn{n+1}$, we get
$$\Delta_C^{AB}\,m_{BD}^E\;\langle \chi_i\,,\,e_A\rangle\;
[\rdn{n}(e^D)]_{ij}=
\Delta_C^{EA}\;\langle \chi_j\,,\,e_A \rangle \ .
\eqno(16)$$
Let us saturate (16) with $e_E\tens e^C$,
use the properties of the Hopf algebras and take into account the
expression of $R_{ij}$ in terms of the double representation as given in
Theorem (6). It is then straightforward to derive
$$(\uno\tens \chi_i)\,T\,(S^{-1}(R_{ji})\tens\uno)=T\,(\uno\tens
\chi_j)\ .
$$
Multiplying to the left by $T^{-1}$ and to the right by
$R_{kj}\tens\uno$, we finally obtain
$$T^{-1}\,(\uno\tens\chi_{k})\,T=R_{ki}\tens\chi_i\ .
$$

This expression is clearly equivalent to the ad-invariance
property of $\chi_i$, {\it i.e.}
$$\ad_X\,\chi_i=\langle X\,,\,R_{ik}\rangle\,\chi_k\ , \quad\quad
X\in \uq\ .\eqno(17)
$$
Therefore the elements
$\chi_i$ verifying (15) and (17) linearly generate a quantum Lie algebra
and define thus a bicovariant differential calculus.

The converse part of this theorem is easily obtained proceeding in the
reverse direction. \fidi
\smallskip
\rmks 18 $(i)$ Saturating (16) with $e_E$ we get
$$(a\star\chi_i)=(\chi_j\star a)\,R_{ij}\ ,\quad\quad a\in \fq\,.
$$
Formally this is the same rule for passing from the left-invariant to
the right-invariant vector fields in Lie group theory.
\smallskip
$\phantom{i}(ii)$ The classical case is obtained by observing that
the double of a Lie group $G$ is $\dop={\cal U}(T^*G)$, where
the brackets between elements of the Lie algebra ${\rm Lie}\,G$ and its
dual $({\rm Lie}\,G)^*$ are
$$[\,\theta^i,X_j\,] = f^i_{jk}\,\theta^k\,,~~~~~~~X_j\in{\rm Lie}\,G\,,
{}~\theta^i\in({\rm Lie}\,G)^*\,.$$
Here $f^i_{jk}$ are the algebra structure constants. Using the
coadjoint representation for ${\rm Lie}\,G$ and the trivial representation
in the same dimension for $({\rm Lie}\,G)^*$,
it is easy to see that the conditions required
in Theorem (13) are satisfied and the classical differential calculus
is easily deduced.
\smallskip
$(iii)$ Representing $\sigma\circ {\cal R}^{-1}$ with $\rdn{n+1}$ we
obtain the matrix
$$\Lambda_{cd}^{ab}=[\rd(e_A)]_{bd}\,[\rd(\sd^{-1}(e^A)]_{ac}\ ,\quad\quad
{\rm with}~~~\ \ a,b,c,d\,=\,0,1,\cdots,\,n
$$ whose nonzero entries are $\Lambda_{k\ell}^{ij}$,
$\Lambda_{k\ell}^{j0}=\langle\chi_\ell\,,\,R_{kj}\rangle$, and
$\Lambda_{b0}^{a0}=\Lambda_{0b}^{0a}=\delta_{ab}$, where
$i,j,k,\ell=1,\cdots,n$.
We recover the structure of the quasi-triangular quantum Lie algebras
defined by Bernard in [7]. We observe that, also in this case, the quasi
triangularity is implied by the connection with the Drinfeld double. \fidi

\bigskip
{\bf 6.}
Let us give the construction of a four dimensional
differential calculus for $E_q(2)$. We emphasize that this
quantum group is not quasi triangular nor is its Lie-Poisson counterpart
coboundary. A differential calculus on a different deformation of
the Euclidean group, non quasi-triangular too, has been obtained in
[8] studying directly the ad-invariant right ideals.

The double of $E_q(2)$ is generated by three elements of the quantum
enveloping algebra $J,b_+,b_-$
and by the corresponding quantized canonical coordinates of the second
kind  $\pi,\pi_+,\pi_-$.
Their duality relationships read
$$\langle J,\pi\rangle =\langle b_+,\pi_+\rangle =
\langle b_-,\pi_-\rangle = 1\,.$$
The algebraic relations for the double are
$$[J,b_+]=b_+\,,~~~~~~~[J,b_-]=-b_-\,,~~~~~~~
[b_+,b_-]=0\,,$$
$$[\pi,\pi_+]=-z\pi_+\,,~~~~~~~[\pi,\pi_-]=-z\pi_-\,,~~~~~~~
[\pi_-,\pi_+]=0\,,$$
together with
$$[b_-,\pi_-]= e^{-\pi}- e^{-zJ}\,,~~~~~~~
[b_-,\pi]=-zb_-\,,~~~~~~~
b_-\pi_+ - e^ z\pi_+ b_-=0\,,$$
$$[J,\pi_-]=\pi_-\,,~~~~~~~[J,\pi]=0\,,~~~~~~~
[J,\pi_+]=-\pi_+\,,$$
$$[b_+,\pi_+]=- e^{-\pi}+ e^{zJ}\,,~~~~~~~
[b_+,\pi]=-zb_+\,,~~~~~~~
b_+\pi_- - e^ z\pi_- b_+=0\,.$$
\smallskip
\noindent The coalgebra of the double is as follows
$$\dd b_+=b_+\tens 1+ e^{zJ}\tens b_+\,,~~~~~~~
\dd b_-=b_-\tens  e^{-zJ}+1 \tens b_-\,,~~~~~~~
\dd J=J\tens 1+1\tens J\,,$$
$$\dd \pi_+=\pi_+\tens  e^{-\pi}+1 \tens \pi_+\,,~~~~~~~
\dd \pi_-=\pi_-\tens 1+ e^{-\pi}\tens \pi_-\,,~~~~~~~
\dd \pi=\pi\tens 1+1\tens \pi\,,$$
while the antipode reads
$$\sd (b_+)=- e^{-zJ}b_+\,,~~~~~\sd (b_-)=-b_-  e^{zJ}\,,~~~~~~~
\sd (J)=-J\,,$$
$$\sd (\pi_+)=-\pi_+ e^{\pi}\,,~~~~~~~\sd (\pi_-)=-  e^{\pi}\pi_-\,,
{}~~~~~~~~~\sd (\pi)=-\pi\,.$$
\smallskip
Bicovariant bimodules are obtained from
representations of the double in dimensions two and three, however
a bicovariant differential calculus is found in dimension four. The
appropriate representation of the double is specified as follows:
$$\rdn{4+1}(J)=-e_{22}+e_{33}\,,~~~~~~~
\rdn{4+1}(b_+)= e^{-3z/4}e_{12}+ e^{z/4}e_{34}\,,$$
$$\rdn{4+1}(b_-)= e^{z/4}e_{13}+ e^{5z/4}e_{24}\,,$$
$$\rdn{4+1}(\pi)=(-z/\kappa)e_{04}+z(e_{11}-e_{44})\,,~~~~~~~
\rdn{4+1}(\pi_+)= e^{z/2}e_{03}+ e^{z/2}\kappa (e_{21}+e_{43})\,,$$
$$\rdn{4+1}(\pi_-)=- e^{-z/2}e_{02}- e^{-z/2}\kappa (e_{31}+e_{42})\,,$$
where $\kappa=2 e^{-z/4}{\rm sh}(z/2)$ and $e_{ij}$ are the usual
matrices with unity in the $(ij)$ entry.

The invariant vector fields $\chi_i$ turn out to be
$$
\chi_1=-\kappa b_-b_+\,,~~~~~
\chi_2=- e^{-z/2}b_-\,,~~~~~\chi_3= e^{z/2} e^{-zJ}b_+\,,~~~~~
\chi_4=\kappa^{-1}( e^{-zJ}-1)\,.$$

Using the universal $T$-matrix
$$T= e^{\pi_-\tens b_-}_z\, e^{\pi\tens J}\, e^{\pi_+\tens b_+}_{-z}\,,$$
from Theorem (6) we get
\medskip
$$(f_{ij})=\left(\matrix{
 e^{zJ} & 0 & 0 & 0\cr
\kappa e^{z/2}b_+ & 1 & 0 & 0\cr
-\kappa e^{-z/2}b_- e^{zJ}  & 0 & 1 & 0\cr
-\kappa^2 b_- b_+ & -\kappa e^{-z/2}b_- & \kappa e^{z/2} e^{-zJ}b_+
&  e^{-zJ}\cr}\right)$$
\medskip
\noindent and
\medskip
$$(R_{ij})=\left(\matrix{
1 & 0 & 0 & 0\cr
- e^{z/4}\bar n & \bar v & 0 & 0\cr
- e^{z/4}n  & 0 & v & 0\cr
 e^{z/2}n\bar n & - e^{z/4}n\bar v & - e^{z/4}v\bar n & 1\cr}
\right)\,,$$
\medskip
\noindent
where $v= e^{-\pi}$, $n=\pi_-$, $\bar n= e^{\pi}\pi_+$ generate
$\fq(E_q(2))$.

\bigskip
{\bf 7.}
In this section we present some cohomological features
of the previous construction of differential calculi. The extension
of the representation as described in Theorem (13) is indeed connected
with a Hochschild cohomology that takes values in the bimodule of
invariant forms.

The representation space $\gin$ of a $n$-dimensional representation $\rd$
can be given a $\dop$-bimodule structure as follows:
$$\alpha\cdot v=\epsilon(\alpha)\, v\,,~~~~~~~v\cdot\alpha=[\rd(\alpha)]^t
\,v$$
with $v\in\gin$, $\alpha\in\dop$.
The structure of $\fq$-bicovariant bimodule on $\Gamma=\fq\tens \gin$
given in (2) is recovered as
$$a\cdot (b\tens v)=a b\tens v\,,\qquad  (b\tens v)\cdot
a=\sum\limits_{(a)}b a_{(1)}\tens v\cdot a_{(2)}\,.$$
We call $\ckd$ the set of $k$-cochains on $\dop$, namely
the $k$-multilinear mappings $\varphi$ from $\dop^k$ to $\gin$, with
$C^0(\dop,\gin)= \gin$.
We then define the coboundary operator $\delta:\ckd\longrightarrow\ck1$ as
$$\eqalign{
(\delta\varphi)\,&(\alpha_1,\,\alpha_2,\,\dots\,,\,\alpha_{k+1})=
\alpha_1\cdot\varphi(\alpha_2,\,\dots\,,\,\alpha_{k+1})+\cr
{}&\sum\limits_{i=1}^{k}\,(-1)^i\,
\varphi(\alpha_1,\dots\,,\,\alpha_i\,\alpha_{i+1},\dots\,,\,
\alpha_{k+1}) +
(-1)^{k+1}\,
\varphi(\alpha_1,\dots\,,\,\alpha_k)\cdot\alpha_{k+1}\,.\cr}
$$
It is a standard fact that $\delta^2=0$. Hence Hochschild cocycles
$Z^k(\dop,\gin)$, coboundaries $B^k(\dop,\gin)$ and
cohomology groups $H^k(\dop,\gin)$ are defined as usual.
\smallskip
Using the explicit expression for 1-cocycles,
$$\delta\varphi(\alpha_1,\alpha_2)=\epsilon(\alpha_1)\,\varphi(\alpha_2)-
\varphi(\alpha_1\,\alpha_2)+[\rd(\alpha_2)]^t\,\varphi(\alpha_1)=0\,,$$
the statement of Theorem (13) can be easily cast into the following form.
\smallskip
\pro{19} {\it Bicovariant differential calculi are in one to one
correspondence with $1$-cocycles $\varphi\in Z^1(\dop,\gin)$ satisfying the
additional condition
$\varphi(X)=0$ for any $X\in\uq$.} \fidi
\smallskip
\rmk{20} ({\it Universal calculus}.)
This approach allows us to describe the universal calculus.
It is useful to write the relations (8) of the double in the intrinsic
form:
$$\eqalign{
X\,a &= \sum\limits_{(a)\,(X)}\,a_{(2)}\,X_{(2)}\,\,\,\langle X_{(1)}\,,\,
a_{(3)}\rangle\,
\,\langle X_{(3)}\,,\,S^{-1}(a_{(1)})\rangle\,,\cr
a\,X &= \sum\limits_{(a)\,(X)}\,X_{(2)}\,a_{(2)}\,\,\,\langle X_{(1)}\,,\,
S^{-1}(a_{(3)})\rangle\,\langle X_{(3)}\,,\,a_{(1)}\rangle\,,\cr}\eqno(21)
$$
($a\in \fq$, $X\in \uq$).
With a direct calculation using (21) it can be seen that
$\,\,{\rm ker}\, \epsilon \subseteq \fq$ is a
$\dop$-bimodule according to
$$\eqalign{
a\cdot h = \epsilon(a)\, h&\,, \qquad h\cdot a = h\, a\,,\cr
X\cdot h = \epsilon(X)\, h&\,, \qquad
h\cdot X = \Ad_{\sd{(X)}}(h)\,, \cr
}$$
where $h\in {\rm ker}\, \epsilon$ and
$\Ad_X(a)=(1\otimes X)\,\ad^*a$.
It turns out that, defining $\fic\in C^1(\dop,{\rm ker}\,\epsilon)$ as
$$\fic(aX)=\Ad_{\sd(X)}(a)-\epsilon(X)\,\epsilon(a)\,,$$
we have
$$\delta\,\fic = 0\,,~~~~~{\rm and}~~~~~\fic(X) = 0\,.$$
Therefore $\fic$ is a 1-cocycle defining a differential calculus. In order
to gain further insight into the result, we recall [1] that the map
$$a\tens b\mapsto r(a\tens b) = (a\tens 1)\,\Delta b :\fq\tens\fq
\rightarrow\fq\tens\fq$$
establishes a bimodule isomorphism of $\fq^2 := {\rm ker}\, m$ with
$\fq\tens{\rm ker}\,\epsilon$, the former with the standard $\fq$-bimodule
structure, the latter with the following one: for $x=\sum\limits_k a_k
\tens h_k$,
$$a\cdot x = \sum\limits_k (a\,a_k)\tens h_k\,,~~~~~~~
x\cdot a= \sum\limits_k a_k\tens h_k\,\,\Delta a\,,$$
($a\in\fq$, $x\in\fq\tens{\rm ker}\,\epsilon$).
Using the cocycle $\fic$ we find a differential
$$D'\,a = \sum\limits_{(a)}\,a_{(1)}\tens\fic(a_{(2)}) = \Delta\,a -
a\tens 1\,.$$
We then see that $D' = r\circ D$, where $D a = 1\tens a - a\tens 1$ is the
differential of the universal calculus. \fidi
\smallskip
Let us show that differential calculi can be specified
in terms of the Hochschild cohomology of the algebra of functions,
obtained by an obvious restriction of the one defined for the double.
The additional
condition necessary to define a differential calculus results into
the invariance of the cocycles under the action of $\uq$ defined on
the cochains as
$$(\psi\bullet X)\,(a_1\,,\,\dots\,,\,a_k) =
\sum\limits_{(X)}\,[\rd(X_{(k+1)})]^t\,
\psi(\Ad_{X_{(k)}}\,a_1\,,\,\dots\,,\Ad_{X_{(1)}}\,a_k)\,,
\eqno(22)$$
with $\psi\in C^k(\fq,\gin)$.
The invariance under this action is, as usual,
$\psi\bullet X=\epsilon(X)\,\psi$. Denote by ${\widetilde C}^0(\fq,\gin)$
the invariant 0-cochains and by ${\widetilde Z}^1(\fq,\gin)$
the invariant 1-cocycles.
\smallskip
\pro{23} {\it $\phantom{i}(i)$ There is a one to one correspondence
between differential calculi and invariant $1$-cocycles
$\psi\in{\widetilde Z}^1(\fq,\gin)$.

$(ii)$ The coboundary operator $\delta$ maps
${\widetilde C}^0(\fq,\gin)$ into ${\widetilde Z}^1(\fq,\gin)$, so that
each invariant $0$-cochain defines a coboundary differential calculus.}
\smallskip
\dim $(i)$ Let $\varphi\in Z^1(\dop,\gin)$ be a 1-cocycle that
determines a differential calculus. According to the previous definitions,
we have
$$\varphi(X\,a)=\epsilon(X)\,\varphi(a)\,,\quad \varphi(a\,X)=[\rd(X)]^t\,
\varphi(a)\,,\qquad a\in \fq\,,\ X\in \uq\,.$$
The second relation, by use of the second of (21), becomes
$\varphi(\Ad_{\sd(X)}(a))= [\rd(X)]^t\,\varphi(a)\,.$
Define $\psi$ to be the restriction of $\varphi$ to $\fq$.
Then
$$\sum\limits_{(X)}\,[\rd(X_{(2)})]^t
\,\psi(\Ad_{X_{(1)}}(a))  = \epsilon(X)\,\psi(a)\,.$$
Moreover $\delta\psi=0$ as a consequence of $\delta\varphi=0$, which
holds by assumption.

Conversely, let $\psi$ be an invariant 1-cocycle. Define $\varphi =
\psi\circ\fic$, where $\fic$ is defined in (20). It is easily seen that
$\varphi(X\,a)=\epsilon(X)\,\varphi(a)$ and $\varphi(X)=0$ for any
$X\in\uq$. If $\psi\bullet X=\epsilon(X)\psi$ then
$$\eqalign{
\psi(\Ad_X(a))&=\sum\limits_{(X)}[\rd(X_{(2)}\,S(X_{(3)}))]^t
\psi(\Ad_{X_{(1)}}(a))\cr
{}&=\sum\limits_{(X)}[\rd(S(X_{(2)}))]^t\epsilon(X_{(1)})
\psi(a)=[\rd(S(X))]^t\,\psi(a)\,.}
$$
We thus have
$$\varphi(a\,X)=\psi(\Ad_{\sd(X)}(a))=[\rd(X)]^t\,\varphi(a)\,,$$
so that $\varphi$ is a 1-cocycle with $\varphi(X)=0$ and defines a
differential calculus.
\smallskip
$(ii)$  Let $\gamma\in {\widetilde C}^0(\fq,\gin)$, namely
$[\rd(X)]^t\,\gamma=\epsilon(X)\,\gamma$ for any $X\in\uq$. We want
to prove the invariance of $\delta\gamma$, {\it i.e.} $\delta\gamma
\bullet X = \epsilon(X)\,\delta\gamma$.
This is equivalent to showing
$$[\rd(X)]^t\,(\delta\gamma)(a)=(\delta\gamma)(\Ad_{\sd(X)}(a))\,.$$
Using the explicit expression for $\delta\gamma$, the previous relation
reads
$$[\rd(aX)]^t\,\gamma=
[\rd(\Ad_{\sd(X)}(a))]^t\,\gamma\,,$$
which holds as a consequence of the second relation of (21). \fidi
\smallskip
\cor{24} {\it Let $\gamma\in {\widetilde C}^0(\fq,\gin)$, $a\in\fq$.
The cochain $\gamma$ is left and right invariant and its coboundary
$\delta\gamma$  induces a differential}
$$da=a\cdot (1\tens \gamma)-(1\tens \gamma)\cdot a\,.$$
\smallskip
\dim The left-invariance of $\gamma$ is by definition. Then, from (5)
and due to the fact that $\gamma\in {\widetilde C}^0(\fq,\gin)$, we have
$(1\tens X)\,\gd\gamma = [\rd(X)]^t\gamma = \epsilon(X)\gamma$.
Hence $\gd\gamma = \gamma\tens 1$ and the right-invariance is proved.

For the second statement we have that
$(\delta \gamma)(a)=(\epsilon(a)-[\rd(a)]^t)\,\gamma$ and
the result comes from $da=\sum_{(a)}\,a_{(1)}
\tens(\delta\gamma)(a_{(2)})\,.$ \fidi

\smallskip
\rmks {25} We shall conclude by showing that
most of the known results on differential calculi can be
organized in a coherent way. We shall also to point out a difference
of the behaviour of the usual differential calculus  on Lie groups.
\smallskip
$\phantom{iiv}(i)$
({\it Quantum groups of type A,B,C,D.}) All the  differential calculi
for the quantization of the simple Lie algebras of the series
$A,\,B,\,C,\,D$ have been classified in [3]. They have been proved
to be {\it inner}, in the sense that there exists a left and right
invariant form $\omega$ such that
$da=a \omega -\omega a$. From Corollary (24), we see that all
those differential calculi are coboundary and determined by
$\delta\omega$.
\smallskip
$\phantom{iv}(ii)$
({\it Quantum group $E_q(2)$.})
The differential calculus described
in the previous section is generated by the coboundary $-{\kappa}^{-1}\,
(\delta \omega_4)$. It is also interesting to observe that this quantum
group has a limiting Lie-Poisson structure that is {\it not} coboundary.
\smallskip
$\phantom{v}(iii)$
({\it Finite groups.}) The approach described in this paper can be used
to determine bicovariant differential calculi on finite groups using the
representations of the double [9]. It is immediate to recover the results
presented in [10]. The differential calculi on a finite group $G$ are in
one to one correspondence with the set $\{{\cal C}\}$ of its
conjugacy classes. The space of
invariant forms is the linear space $\langle\omega_g\rangle\,,\
g\in{\cal C}$ which carries the representation of the double
$\dop(G)$
$$\rd(h)\,\omega_g=\omega_{hgh^{-1}}\,, \qquad \rd(a)\,\omega_g=\langle
g\,,\,a\rangle \, \omega_g\,,$$
where $h\in G$, $a\in\fq(G)$ and the dual pairing is defined by
$\langle g\,,\,a\rangle = a(g)$.
The invariant 1-cocycle $\psi$ that defines the differential
calculus has components  $[\psi(a)]_g=\epsilon(a)-\langle
g\,,\, a\rangle\,$, with $a\in \fq(G)$, $g\in {\cal C}$.
Also in this case $\psi$ is a coboundary, namely
$\psi=\delta(\sum_{{\cal C}} \omega_{g})$. \fidi
\smallskip
$\phantom{ii}(iv)$
({\it The classical differential calculus on Lie groups.})
{}From the Remark (18) the
representation of $\fq$ corresponding to the differential calculus is
the trivial one, namely $\rd (a)=\epsilon(a)\,,\forall a\in \fq$.
We then have that $(\delta\gamma)(a)=(\epsilon(a)-[\rd(a)]^t)\,\gamma=0$.
Hence the coboundaries are zero and
the classical differential calculus is associated with a
nontrivial 1-cocycle. Therefore, at difference with
cases $(i)-(iii)$, the differential calculus for Lie groups
is not inner.

\bigskip
\bigskip
\baselineskip= 10 pt
{\refbf Acknowledgement.} {\abs We are indebted to G. Vezzosi for
useful discussions and to the referee for constructive criticism.
After this work was completed S. Zakrzewski informed us that
the content of Theorem (6) was present in the paper ``{\refit Quantum
deformation of Lorentz group}'' by P. Podl\'es and S.L. Woronowicz,
{\refit Comm. Math. Phys.}, {\refbf 130} (1990) 381.
}

\bigskip
\bigskip
\bigskip

\centerline{{\bf References.}}

\bigskip
\baselineskip= 10 pt
{\abs
\ii 1 S.L. Woronowicz, {\refit Commun. Math. Phys.}, {\refbf 122}
(1989) 125.
\smallskip
\ii 2 B. Jurco, {\refit Lett. Math. Phys.}, {\refbf 22} (1991) 177.
\smallskip
\ii 3 K. Schm\"udgen and A. Sch\"uler, ``Classification of
bicovariant differential calculi on quantum groups of type A, B, C and D'',
Leipzig preprint KMU-NTZ 94-1.
\smallskip
\ii 4 V.G. Drinfeld, {\refit Proc. Int. Congr. Math., Berkeley} (1986)
798.
\smallskip
\ii 5 L.L. Vaksman and L.I. Korogodski, Sov. Math. Dokl. {\bf 39},
173 (1989).
\smallskip
\ii 6 N.Yu Reshetikhin and M.A. Semenov-Tian-Shansky, {\refit
Journal Geom. Phys.}, {\refbf 5} (1988) 533.
\smallskip
\ii 7 D. Bernard, {\refit Prog. Theor. Phys. Supp.},
{\refbf 102} (1990) 49.
\smallskip
\ii 8  S. Giller, C. Gonera, P. Kosinski, P. Maslanka,
``Differential Calculus on deformed E(2) Group'', q-alg 9501011 (1995).
\smallskip
\ii 9 R. Dijkgraaf, V. Pasquier and P. Roche, `` Quasi Hopf
algebras, group cohomology and orbifold models'', in {\refit Integrable
Systems and Quantum Groups}, eds. M. Carfora, M. Martellini and A.
Marzuoli (World Scientific, Singapore, 1992).
\smallskip
\jj {10} K. Bresser, F. M\"uller-Hoissen, A. Dimakis and A. Sitarz,
``Noncommutative geometry of finite groups'',
 GOET-TP 95/95, q-alg/9509004 (1995).
\smallskip
}

\bye